\documentclass[a4paper,10pt]{article}

\def \Z {{\mathbf {Z}}}
\def \R {{\mathbf {R}}}

\title{On isomorphism of tensor powers of ergodic flows}
\author{V.V. Ryzhikov}
\date{}
\begin{document}

\maketitle

\begin{abstract}{\large
The following question due to Thouvenot is well-known in ergodic theory. 
Let $S$ and $T$ be automorphisms of a probability space and 
$ S \otimes S $ be isomorphic to $T \otimes T $. Will $ S $ and $ T $ be isomorphic? 
Our note contains a simple answer to this question and a generalization
of Kulaga's result on the corresponding isomorphism within a class of  flows. 
We show that the isomorphism of weakly mixing flows $ S_t \otimes S_t $ and 
$ T_t \otimes T_t $ implies the isomorphism of the flows $S_t$ and $T_t$, 
if one of them has an integral weak limit.}
\end{abstract}

\Large

\section{Introduction}
The work is devoted to the problem of the isomorphism of two dynamical systems 
under the condition that  their tensor powers are isomorphic.
 By isomorphisms we mind  a conjugation by  an invertible measure-preserving transformation 
(an automorphism).
 In [1] it was proved that for generic transformations, 
the isomorphism of their tensor powers implies the  isomorphism of the transformations themselves. 
From the results of [3,4] one can deduce a similar result for 
typical actions of multi-parameter
 flows.

We generalize one of the facts from [5] on  the  isomorphism of flows and give 
a negative response 
to the mentioned question in the situation 
where  systems are not required
 to both be ergodic.
  We first start with these examles.

\bf Counterexamples. \rm Let's consider an irrational  shift of the 
circle $X=\R/\Z$. 
 It preserves the Lebesgue measure $\mu$ and it is ergodic.
 The ergodicity means that any  $R$-invariant  measurable function on $X$ 
has to be constant.

The torus $X\times X$ is stratified into   $R\otimes R$-invariant circles, where 
 $R\otimes R(x,y)=(Rx,Ry)$. 
The transformation $R\otimes R$ acts on each circle 
as the original rotation  $R$. 

An isomorphism of the transformations $R\otimes R$ and $I\otimes R$, where 
$I$ denotes the identity transformation, 
is given by  the map 

$$(x,y)\to \left(x-y, \frac{x+y}{2}\right).$$ 

Setting
$S=I\otimes R$ and $T=R$, 
we get:
 $S\otimes S$ and $T\otimes T$ are isomorphic. Indeed, the transformation 
$I\otimes R\otimes I\otimes R$ is isomorphic to  
$I\otimes I\otimes I\otimes R$, 
 hence, it isomorphic to the transformation $I\otimes R$. However the non-ergodic transformation 
$S$ is 
not isomorphic to 
ergodic $T$.

Similar counterexamples could be obtained with arbitrary ergodic transformation $R$ 
with purely  discrete 
spectrum. 
In the case of flows, let us consider  
$R_t$, an  
 ergodic torus winding, and put $T_t=R_t$ and $S_t=I\otimes R_t$.
  All tensor powers of the flows
$T_t$, $S_t$ are isomorphic between themselves, but the flows $T_t$, $S_t$ are not isomorphic, 
since the first flow is ergodic, and the second is not.

Thouvenot's question remains open 
in the class of weakly mixing
 dynamical systems, in particular, for the flows with continuous
 spectrum.
  Counterexamples within this  class, if they exist, must have, in our opinion, 
unusual properties. The following question is also of interest: 
\it for which groups
 $G$ there are no  such counterexamples  among the measure-preserving G-actions? \rm
 
\bf Main result. \rm In  [5], in particular, it is proved that
the isomorphism of flows  $S_t\otimes S_t$ и $T_t\otimes T_t$ implies the isomorphism of the flows 
 $S_t$ and $T_t$, if the flow $T_t$ has a weak limit in the form
$\int_{\bf R} T_a d\nu(a),$
where $\nu$ is a continuous mesure on $\bf R$ with analytical Fourier transform.
 We shall show that this last condition is not necessary.

\vspace{4mm}

\bf Theorem. \it 
The isomorphism of flows  $S_t\otimes S_t$ and $T_t\otimes T_t$ implies the isomorphism of the flows 
 $S_t$ and $T_t$, if the flow  $T_t$ has a weak limit in the form
$$\int_{\bf R} T_a d\nu(a),$$
where $\nu$ is a continuous mesure on $\bf R$. \rm

\vspace{4mm}
Remark. The theorem is true for any non-Dirac measure $ \nu $.
The presence of pointwise components of the measure simplifies the proof (see [2]).
 We confine ourselves to the most interesting case.

\section{ Proof of Theorem} 

Let $ \Phi $ denote the isomorphism of the  measure-preserving flows $ S_t \otimes S_t $ and $ T_t \otimes T_t $.
Subsequently, the transformations and the operators on $L_2$ corresponding to them 
are denoted identically.
Instead of equations of the form
$$ S_t \otimes S_t = \Phi (T_t \otimes T_t) \Phi^{- 1} $$
below we will write
$$ S_t \otimes S_t =_\Phi T_t \otimes T_t. $$
For some sequence $ {t_i} $, we have
$$ T_{t_i} \ \to \ \int_{\bf R} T_a d \nu (a) $$
(here and below we consider the weak operator convergence).
For some Markov operator $ Q $ we have
$$ S_{t_i} \ \to \ Q $$
and the equality
$$ Q \otimes Q \ =_\Phi \ \int \int T_a \otimes T_b d \nu (a) d \nu (b). $$
Markov operator $ Q $ commutes with  the flow $S_t $. 
Recall that  Markov operators
preserve the non-negativity of functions and send the constants to themselves. 
Let us  consider 
$ \rho $, the measure on $ X \times X $, defined by the relation
$$ \rho (A \times B) = (Q \chi_A, \chi_B) $$
for all measurable sets $ A, B $.
It is invariant with respect to the flow
 $ S_t \otimes S_t $ and has the standard marginals: its  projections onto the factors in
$ X \times X $ are  $\mu$. In ergodic theory such measures are called self-joinings. 

The measure $ \rho $ decomposes into ergodic with respect to $ S_t \otimes S_t $
components $ \rho_c, c\in C,$  
$$\rho=\int_C \rho_c\sigma(c).$$ 
  Let us show the ergodicity (almost surely) of the measures $ \rho_c \times \rho_{c'} $
with respect to the transformation $ S_t \otimes S_t \otimes S_t \otimes S_t $.
Recall the well-known fact of the spectral theory of ergodic transformations and flows:
the tensor product of  two systems is not ergodic iff these systems  have the same 
eigenvalue different from 1.
This eigenvalue is inherited by ergodic components of this tensor product.

Let us consider the self-joining of $ \eta $, corresponding to the operator
$$ \int_R \int_R T_a \otimes T_b \ d \nu (a) d \nu (b). $$
The ergodic components of $ \eta $ are sitting on the graphs
of the transformations
$ T_a \otimes T_b $.
The flow $ (T_t \otimes T_t \otimes T_t \otimes T_t, \eta) $
is isomorphic to the flow $ (S_t \otimes S_t \otimes S_t \otimes S_t, \rho \times \rho) $.
All ergodic components of the flow $ (T_t \otimes T_t \otimes T_t \otimes T_t, \eta) $ are isomorphic
to the weakly mixing flow $ (T_t \otimes T_t, \mu \times \mu) $.
Consequently,
for almost all  $c,{c '}$ with respect to $\nu\times\nu$ the 
dynamical system
$ (S_t \otimes S_t \otimes S_t \otimes S_t, \rho_c \times \rho_{c '}) $
is isomorphic via $\Phi$  to the flow $ (T_t \otimes T_t, \mu \times \mu) $. Indeed, if the measures
$ \rho_c \times \rho_{c '} $ is non-ergodic, then the flow
$ (T_t \otimes T_t, \mu \times \mu) $ possesses a  non-constant proper function, but it is not true.

Let $ Q_{c} $ denote the Markov operators, corresponding to the self-joining  $ \rho_{c} $.
We have
$$Q=\int_C Q_c d\sigma(c),$$
$$\int_C \int_C Q_c\otimes Q_{c'}d\sigma(c)d\sigma(c')=_\Phi \int_R \int_R T_a\otimes T_b \ d\nu(a) d\nu(b).$$ 
The isomorphism $ \Phi $ maps extreme points to extreme points, so, it maps 
the Markov operators $ Q_{c} \otimes Q_{c '} $ to
the operators $ T_a \otimes T_b $.
The isometry maps  automorphisms  to  automorphisms, so $ Q_{c} $ is an automorphism
commuting with the flow  $ S_t $. For some functions ${c(a, b)}$ and ${c'(a, b)}$ we obtain  the equality
$$ Q_{c(a, b)} \otimes Q_{c'(a, b)} =_\Phi T_a \otimes T_b $$
for almost all $ a, b $ with respect to the measure  $ \nu $.
From
$ S_t \otimes S_t =_\Phi T_t \otimes T_t, $
we get now
$$ S_{- b} Q_{c (a, b)} \otimes S_{- b} Q_{c'(a, b)} =_\Phi T_{a-b} \otimes I. \eqno (\ast ) $$
This equality holds for some different numbers $ a, b$.  We pick them.

\bf Case 1. \rm If one of the factors in the product $ S_{- b} Q_{c} \otimes S_{- b} Q_{c '} $
is the identity operator,
and the other is ergodic, then  Theorem is proved. Indeed, in this case the isomorphism $ \Phi $
 maps the coordinate
algebra, corresponding to the identity operator  from  $ T_{a-b} \otimes I $,
 to one of the  coordinate algebras. It remains to show that other cases cannot be realized.

\bf Case 2. \rm 
Let  $S_{-b}Q_{c}$ and $S_{-b}Q_{c'}$ both be  ergodic. 
They will have a continuous spectrum. 
Indeed, 
the automorphism $T_{a-b}$  has continuous spectrum, therefore, the operator
 $T_{a-b}\otimes I$ has no eigenvalues,
 except  1. So the product  
 $S_{-b}Q_{c}\otimes S_{-b}Q_{b'}$ has continuous spectrum,
 it cannot be isomorphic
to $T_{a-b} \otimes I.$

\bf Case 3. \rm Let  the  automorphisms of $S_{-b}Q_{c}$ 
and $S_{-b}Q_{c'}$ both 
be not ergodic and each of them is not the identity automorphism.

Following Rokhlin, we represent them in the form of skew products over identity transformations.

From  $(\ast)$ we see  that the automorphism $T_{a-b}$ can be represented as $T'\otimes T"$,
where $T'$ and $T''$ 
are isomorphic to the ergodic components of the automorphisms  
$S_{-b}Q_{c}$ and $S_{-b}Q_{c'}$, respectively.
  We will show that this case is impossible.

Put $d=b-a\neq 0$.
 We have
$$T_{t_i}\to\int_{\bf R} T_a d\nu(a).$$  We find  some integer sequence 
$k_i$ such that 
 $t_i= k_id+r_i$, $r_i\to r\in [0,d)$,
and
$$T_{k_id}\ \to \ J=\int_{\bf R} T_u d\nu'(u),$$
where $\nu'$ is $r$-shift of the measure  $\nu$. 
The limit operator $J$  (as well as all  operators $T_{k_id}$) commutes
with 
the orthogonal projection  $P$ 
onto the space of $T'$-factor. Let us rewrite the equality

$$P J=JP$$ in the form

$$\int_{\bf R} P T_u d\nu'(u)=\int_{\bf R} T_vP d\nu'(v).$$


The operators $P T_u, T_vP$ are  extreme points in the convex compact set of 
all Markov operators,
 commuting with the automorphism $T_{d}$.
 This  is equivalent 
to the fact that they correspond to 
ergodic self-joinings. The dynamical system associated 
with these self-joinings is isomorphic to the ergodic transformation $T'\otimes T''\otimes T''$.
 The ergodicity  follows from the property of the weak mixing of the factors.

Thus, for almost all $u$ with respect to the measure $\nu'$ 
there 
is a number $v=v(u)$ such that 
$$T_u P = T_vP.$$

Suppose that there exists $u$ such that $u\neq v(u)$. 
Then we have  
$$T_{u-v}P = T_uPT_{-u},$$

$$T_{(u-v)n}P = T_{un}PT_{-un}.\eqno (\ast \ast)$$

Let us represent 
$$un=dm_n+s_n, \ 0\leq s_n\leq d,$$

and let   $$s_{n_k}\to s$$ and 
 $$T_{(u-v)n_k}\to \Theta,$$ 
where $\Theta$ is the orthoprojection onto the constants in $L_2(X,\mu)$.
 For the weakly mixing automorphism $T_{(u-v)}$ we   find easily such a sequence ${n_k}$.
 Using  

$$P=T_{dm}P T_{-dm}$$ 
from  $(\ast \ast)$ we get 

$$\Theta =\Theta P= T_{s}PT_{-s},$$

  $$\Theta = P,$$ but
 the latter contradicts  the  definition of the operator $P$.

Thus, the assumption that $u\neq v(u)$
 leads to a contradiction. So 
$u=v(u)$, and, hence,  the equality 
$$P T_u = T_uP$$ holds for almost all $u$ with respect  to the  measure $\nu'$.
From the continuity of the flow and the the continuity of the  measure $\nu'$  it follows that
the operator  $P$ commutes with all elements of the flow.
  Therefore, the corresponding  
factor-algebra of the automorphism  $T_d$
  is a factor of the flow $T_{t}$ as well.

The weak limits are inherited by the factors of flow,
thus, from  $$T_{t_i}\to\int_{\bf R} T_u d\nu(u)$$  we get

$$\int_{\bf R} (T_c'\otimes T_c'') d\nu(c)=\int_{\bf R} T_a' d\nu(a)\ \otimes \ \int_{\bf R} T_b'' d\nu(b).$$

But this equality is possible only in the case when the measure $\nu$ is concentrated at one point.

Otherwise we  would have for  some different $a,b$
 the equality
$$ T_c'\otimes T_c''= T_a' \otimes T_b'',$$ which is impossible. 
Thus, the case 3 does not occur.

\bf Case 3'. \rm Recall that in the  case 3 the automorphisms 
 $S_{-b}Q_{c}$ and $S_{-b}Q_{c'}$ 
are considered both non-ergodic and  not the identity automorphism.

Let now  $S_{-b}Q_{c'}$ be  ergodic, then  
it can play the role 
of  $T''$,
 and this situation is actually not different from the case 3. 
So, our present case  is impossible due to the same reasons.

\bf Case 4. \rm Let one of the automorphisms $S_{-b}Q_{c}$ and $S_{-b}Q_{c'}$ 
be non-ergodic and the other
 be the identity operator.
  Applying  

$$ S_{b-a}\otimes 
S_{b-a}=_\Phi T_{b-a}\otimes T_{b-a},$$ from $(\ast)$ we get

$$S_{-a}Q_{c}\otimes S_{-a}Q_{c'}=_\Phi I\otimes T_{b-a}.$$

One of the automorphisms $S_{-a}Q_{c}$, $S_{-a}Q_{c'}$ 
is ergodic, and 
the other is not.
 We  are back to the case 1, since  3' is not possible.
The theorem is proved.

\vspace{5mm}

\normalsize
[1]. V. V. Ryzhikov. Genericity of the isomorphism of measure-preserving
transformations under isomorphism of their Cartesian powers. Mat.
Zametki, 59(4):630-632, 1996.

[2]. V. V. Ryzhikov, A. E. Troitskaya. The tensor root of an isomorphism
and weak limits of transformations. Mat. Zametki, 80(4):596-600, 2006.

[3].   S. V. Tikhonov, Embedding lattice actions in flows with multidimensional time,
 Mat. Sb., 197:1, 97-132, 2006. 

[4]. A. E. Troitskaya. On the isomorphism of measure-preserving $Z^2$-actions 
that have isomorphic Cartesian powers. Fundam. Prikl. Mat., 13(8):193-212, 2007.

[5].  J. Kulaga, A Note on the Isomorphism of Cartesian Products of Ergodic Flows, 
J. Dyn. Control Syst., 18:2, 247-267, 2012.
\vspace{4mm}

vryzh@mail.ru
\end{document}